\documentclass[12pt,leqno]{amsart}
\usepackage{amsmath,amssymb,array,latexsym}
\newtheorem{thm}{Theorem}[section]
\newtheorem{lem}[thm]{Lemma}
\newtheorem{prop}[thm]{Proposition}

\theoremstyle{definition}
\newtheorem{defin}[thm]{Definition}
\newtheorem*{remark}{Remark}

\newtheorem*{question}{Question}
\def\hangbox to #1 #2{\vskip1pt\hangindent #1\noindent \hbox to #1{#2}$\!\!$}

\newcommand{\N}{\mathbb N}

\newcommand{\abs}[1]{\lvert#1\rvert}
\newcommand{\norm}[1]{\lVert#1\rVert}

\newcommand{\Bignorm}[1]{\Bigl\lVert#1\Bigr\rVert}
\newcommand{\trinorm}[1]{\lvert\!\lvert\!\lvert#1\rvert\!\rvert\!\rvert}
\renewcommand{\leq}{\leqslant}
\renewcommand{\geq}{\geqslant}\usepackage{amssymb}

\DeclareMathOperator*{\tlim}{\tau-lim}
\DeclareMathOperator{\dist}{dist}
\DeclareMathOperator{\supp}{supp}
\DeclareMathOperator{\lin}{lin}

\author{Thomas Schlumprecht and Vladimir G. Troitsky}
\title{On quasi-affine  transforms of Read's operator}
\thanks{The first author was supported by NSF. Most of the work on the
  paper was done during the {\it Workshop on linear analysis and
    probability} at Texas A$\&$M University, College Station.}

\begin{document}
\begin{abstract}
 We show that C.~J.~Read's example \cite{Read:85,Read:86} of an
  operator $T$ on $\ell_1$ which does not have any non-trivial
  invariant subspaces is not the adjoint of an operator on a predual
  of $\ell_1$. Furthermore, we present a bounded diagonal operator $D$
  such that even though $D^{-1}$ is unbounded but $D^{-1}TD$ is a
  bounded operator with invariant subspaces, and is adjoint to an
  operator on $c_0$.
\end{abstract}
\maketitle

\markboth{TH. SCHLUMPRECHT AND V. TROITSKY}{ON QUASI-AFFINE TRANSFORMS OF READ'S OPERATOR}
                                            
\baselineskip=15pt

\section{Introduction}\label{S0}
In this note we deal with the Invariant Subspace Problem, the problem
of the existence of a closed non-trivial invariant subspaces for a
given bounded operator on a Banach space. The problem was solved in
the positive for certain classes of operators, see
\cite{Radjavi:73,Abramovich:98} for details, however in the
mid-seventies P.~Enflo \cite{Enflo:76,Enflo:87} constructed an example
of a continuous operator on a Banach space with no invariant
subspaces, thus answering the Invariant Subspace Problem for general
Banach spaces in the negative. In \cite{Read:85} C.~J.~Read presented
an example of a bounded operator $T$ on $\ell_1$ with no invariant
subspace.  Recently V.~Lomonosov suggested that every adjoint operator
has an invariant subspace. In the first part of this note we show that
the Read operator $T$ is not an adjoint of any bounded operator
defined on some predual of $\ell_1$.

Suppose that $A$ has a non-trivial invariant (or a hyperinvariant)
subspace, and suppose that $B$ is similar to $A$, that is,
$B=CAC^{-1}$for some invertible operator $C$. Clearly, $B$ also has a
non-trivial invariant (respectively hyperinvariant) subspace.
Moreover, it is known (see~\cite{Radjavi:73}) that if $A$ has a
hyperinvariant subspace and $B$ is quasi-similar to $A$, (that is,
$CA=BC$ and $AD=DB$, where $C$ and $D$ are two bounded one-to-one
operators with dense range), then $B$ also has a hyperinvariant
subspace. To our knowledge it is still unknown whether or not $A$ has
a non-trivial invariant subspace if and only if $B$ has a non-trivial
invariant subspace, assuming $A$ and $B$ are quasi-similar.

Recall (cf.\cite{Sz-NF:68}) that an operator $A$ is said to be a {\em
  a quasi-affine transform of } $B$ if $CA=BC$, for some injective
operator $C$ with dense range.  In the second part of this paper we
construct an injective diagonal operator $D$ on $\ell_1$ such that even
though $D^{-1}$ is unbounded, the operator $S=D^{-1}TD$ ($T$ being
Read's operator) is bounded and has an invariant subspace.  Thus, we
show that a quasi-affine transform of an operator with no non-trivial
invariant subspace might have a non-trivial invariant subspace.
Furthermore, $S$ is the adjoint of a bounded operator on $c_0$.

Although we prove our statement for a specific choice of $D$, it is
true for a much more general choice, and it seems to be true for any
diagonal operator $D$ that $S=D^{-1} T D$ has a non-trivial invariant
subspace, whenever $S$ is an adjoint of an operator on $c_0$.  More
generally, the following question is of interest in view of above
mentioned conjecture by V.~Lomonosov.

\begin{question} 
  Does every quasi-affine transform of Read's operator, which is an
  adjoint of an operator on $c_0$, have a non-trivial invariant
  subspace?
\end{question}

We introduce the following notations. Following \cite{Read:86} we
denote by $F$ the vector space of all eventually vanishing scalar
sequences, and by $(f_i)$ the standard unit vector basis of $F$.  For
an $x=\sum a_i f_i\in F$, we define the {\em support of $x$} to be the
set $\{ i\in\N: a_i\not=0\}$ and denote it by $\supp(x)$. The linear
span of some subset $A$ of a vector space is denoted by $\lin A$.

\section{Read's operator is not adjoint}\label{S:2}

We begin with reminding the reader of the construction of the operator $T$
in~\cite{Read:85,Read:86}. It  depends on a strictly increasing
 sequence  ${\bf d}=(a_1,b_1,a_2,b_2,\ldots)$  of positive integers 
which has to be chosen to be {\em sufficiently rapidly increasing}.
 Also let $a_0=1$, $v_0=0$,
and $v_n=n(a_n+b_n)$ for $n\geq 1$. 

Read's operator $T$ is defined by prescribing the orbit
$(e_i)_{i\geq 0}$ of the first basis element~$f_0$. 

\begin{defin}\label{d:e_i}
  There is a unique sequence $(e_i)_{i=0}^\infty\subset F$ with the
  following properties:
  \begin{enumerate}
    \item[0)] $f_0=e_0$; 
    \item[A)] if integers $r$, $n$, and $i$ satisfy $0<r\leq n$,
      $i\in[0,v_{n-r}]+ra_n$, we have 
          $$f_i=a_{n-r}(e_i-e_{i-ra_n});$$
    \item[B)] if integers $r$, $n$, and $i$ satisfy $1\leq r<n$,
      $i\in(ra_n+v_{n-r},(r+1)a_n)$, (respectively, $1\leq n$,
      $i\in(v_{n-1},a_n))$, then 
      $$f_i=2^{(h-i)/\sqrt{a_n}}e_i,
\text{ where  $h=(r+\tfrac{1}{2})a_n$ (respectively, $h=\tfrac{1}{2}a_n$);}$$ 
    \item[C)] if integers $r$, $n$, and $i$ satisfy $1\leq r\leq n$,
      $i\in[r(a_n+b_n),na_n+rb_n]$, then 
           $$f_i=e_i-b_ne_{i-b_n};$$
    \item[D)] if integers $r$, $n$, and $i$ satisfy $0\leq r<n$,
      $i\in(na_n+rb_n,(r+1)(a_n+b_n))$, then 
    $$f_i=2^{(h-i)/\sqrt{b_n}}e_i,
     \text{ where $h=(r+\tfrac{1}{2})b_n$.}$$
\end{enumerate}
\end{defin}

Indeed, since $f_i=\sum_{j=0}^{i}\lambda_{ij}e_j$ for each $i\geq 0$
and $\lambda_{ii}$ is always nonzero, this linear relation is
invertible. Further,
\begin{displaymath}
  \lin\{e_i\mid i=1,\dots,n\}=
  \lin\{f_i\mid i=1,\dots,n\}\mbox{ for every }n\geq 0.
\end{displaymath}
In particular, all $e_i$ are linearly independent and also span $F$.
Then C.~J.~Read defines $T\colon F\to F$ to be the unique linear map such that
$Te_i=e_{i+1}$. C.~J.~Read proves that $T$ can be extended to a bounded
operator on $\ell_1$ with no invariant subspaces provided ${\bf d}$
increases sufficiently rapidly.

\begin{prop}\label{P:1.2}
  $T$ is not the adjoint of an operator $S:X\to X$ where $X$ is a
  Banach space whose dual is isometric to $\ell_1$.
\end{prop}

\begin{proof}
  Assume that our claim were not true. Then there is a local convex
  topology $\tau$ on $\ell_1$ so that
  \begin{enumerate}
    \item[a)] $\tau$ is weaker than the norm topology of $\ell_1$;
    \item[b)] B$(\ell_1)$ is sequentially compact with respect to $\tau$;
    \item[c)] If $(x_n)\subset \ell_1$ converges with respect to $\tau$ to $x$, then 
            $\liminf_{n\to\infty}\norm{x_n}\geq\norm{x}$; \label{i:liminf}
    \item[d)] $T$ is continuous with respect to $\tau$.
  \end{enumerate}
   Note that  with respect to any predual $X$ of $\ell_1$ the weak$^*$ topology
   has properties (a) through (d).  
  Let $s\in\N$ be fixed, and
  $n=s,s+1,\ldots$.   Then $f_{(n-s)a_n}=a_s(e_{(n-s)a_n}-e_0)$ by (A)
  above. It follows that
  $T^{v_s+1}f_{(n-s)a_n}=a_s(e_{(n-s)a_n+v_s+1}-e_{v_s+1})$. Further,
  it follows from (B) that $e_{(n-s)a_n+v_s+1}$ equals
  $2^{(1+v_s-\frac{1}{2}a_n)/\sqrt{a_n}}f_{(n-s)a_n+v_s+1}$ and
  converges to zero in norm (and, hence, in $\tau$) as $n\to\infty$.
  Therefore
  \begin{equation}\label{E:1.2.1}
   \tlim_{n\to\infty}T^{v_s+1}f_{(n-s)a_n}=-a_se_{v_s+1}=T^{v_s+1}(-a_se_0).
  \end{equation}
  Notice that $T^{v_s+1}$ is $\tau$-continuous and one-to-one because
  its null space is $T$-invariant. By sequential compactness of $B(\ell_1)$, 
  the sequence $f_{(n-s)a_n}$ must have a
  $\tau$-convergent subsequence, then, by (\ref{E:1.2.1}),  the limit point has to be $-a_s
  e_0$. Since that argument applies to any subsequence, we deduce that
  \begin{equation}\label{E:1.2.2}
    \tlim\limits_{n\to\infty}f_{(n-s)a_n}=-a_s e_0.
  \end{equation}
  Since $\norm{f_{(n-s)a_n}}=1$ for each $n$ and $s$ while
    $\norm{a_se_0}=a_s>1$, this contradicts (\ref{i:liminf}). 
\end{proof}

\begin{remark}\label{R:1.3}
  The statement of the theorem remains valid if we consider an
  equivalent norm on $\ell_1$. Indeed, suppose
  $\frac{1}{K}\trinorm{\cdot}\leq\norm{\cdot}\leq K\trinorm{\cdot}$,
  then $\trinorm{f_{(n-s)a_n}}\leq K$ for each $n$ and $s$, but since
  $\lim_{n\to\infty}a_n=\infty$, we can choose $a_s$ in~(\ref{E:1.2.2})
  so that $\trinorm{a_s e_0}>K$.
\end{remark}

\section{An adjoint operator with invariant subspaces of the form
  $D^{-1}TD$}

Define a sequence of positive reals $(d_i)$ as follows:

\begin{equation}\label{E:2.1}
 d_i=\begin{cases}  
 \frac{1}{r}   &\text{if $ra_m\leq i\leq ra_m+v_{m-r}$ for some $0<r\leq m$}\\ 
 1 &\text{otherwise} 
 \end{cases}
\end{equation}

Let $D$ be the diagonal operator with diagonal $(d_i)$, that is,
$Df_i=d_if_i$ for every $i$. Define $S=D^{-1}TD$. Clearly, $S$ is
defined on $F$. Once we will write $S$ in matrix form it will be clear
that it is bounded on $F$ and, therefore, can be extended to $\ell_1$.
Let $\hat e_i=D^{-1}e_i$, in particular $\hat e_0=e_0$. Then $S\hat
e_i=D^{-1}Te_i=\hat e_{i+1}$, so that the sequence $(\hat e_i)$ is the
orbit of $e_0$ under $S$.

Next, we examine Definition~\ref{d:e_i} to represent the $f_i$'s in terms
of $\hat e_i$'s.

\begin{enumerate}
    \item[$\widehat0$)] $f_0=e_0=\hat e_0$;
    \item[$\widehat{\text{A}}$)] 
    If $i$ satisfies $i\in[0,v_{n-r}]+ra_n$ for some $0<r\leq n$
      then
      $$f_i=d_iD^{-1}f_i=d_iD^{-1}\bigl(a_{n-r}(e_i-e_{i-ra_n})\bigr)=
       \tfrac{a_{n-r}}{r}(\hat e_i-\hat e_{i-ra_n})$$
    \item[$\widehat{\text{B}}$) ] if integers $r$, $n$, and $i$ satisfy $1\leq r<n$,
      $i\in(ra_n+v_{n-r},(r+1)a_n)$, (respectively, $1\leq n$,
      $i\in(v_{n-1},a_n))$, then 
      $$f_i=d_iD^{-1}f_i=2^{(h-i)/\sqrt{a_n}}\hat e_i,\text{ where
      $h=(r+\tfrac{1}{2})a_n$ (respectively, $h=\tfrac{1}{2}a_n$)};$$
    \item[$\widehat{\text{C}}$)] if integers $r$, $n$, and $i$ satisfy $1\leq r\leq n$,
      $i\in[r(a_n+b_n),na_n+rb_n]$, then 
      $$f_i=d_iD^{-1}f_i=\hat e_i-b_n\hat e_{i-b_n};$$
    \item[$\widehat{\text{D}}$)] if integers $r$, $n$, and $i$ satisfy $0\leq r<n$,
      $i\in(na_n+rb_n,(r+1)(a_n+b_n))$, then 
      $$f_i=d_iD^{-1}f_i=2^{(h-i)/\sqrt{b_n}}\hat e_i,
      \text{ where $h=(r+\tfrac{1}{2})b_n$.}$$
\end{enumerate}

We see that it differs from Definition~\ref{d:e_i} only in case ($\widehat{\text{A}}$).
Now we can actually write the matrix of $S$:

\begin{displaymath}
 Sf_i=\begin{cases}
 2^{(1-\frac12a_1)/\sqrt{a_1}}f_1
             &\text{if $i=0$}\\
 f_{i+1}     &\text{if $i\in[0,v_{n-r})+ra_n$},\\
             &\hfil\text{with $r=1,2,\ldots,n$}\\
 f_{i+1}     &\text{if $i\in[r(a_n+b_n),na_n+rb_n)$,}\\
             &\hfil\text{with $r=1,2,\ldots,n$}\\
 2^{1/\sqrt{a_n}}f_{i+1} 
             &\text{if $i\in(ra_n+v_{n-r},(r+1)a_n-1)$,}\\
             &\hfil\text{with $r=1,2,\ldots,n-1$}\\
             &\hfil\text{or $i\in(v_{n-1},a_n-1)$}\\
 2^{1/\sqrt{b_n}}f_{i+1} 
             &\text{if $i\!\in\!(na_n\!+\!rb_n,\!(r\!+\!1)(a_n\!+\!b_n)\!-\!1)$}\\
             &\hfil\text{with $r=0,1,\ldots,n-1$}\\
 \tfrac{a_{n-r}}{r}(\varepsilon_1f_{i+1}\!-\!\varepsilon_2f_{v_{n-r}+1}) 
             &\text{if $i=ra_n+v_{n-r}$,} \\
 \quad\text{where}          
             &\hfil\text{with $r=1,2,\ldots,n$}  \\
 \qquad\varepsilon_2\!=\!2^{(1+v_{n-r}-\frac12a_{n-r+1})/\sqrt{a_{n-r+1}}}&\\
 \qquad\varepsilon_1\!=\!2^{(1+v_{n-r}-\frac12a_{n})/\sqrt{a_n}}   
             &\text{if $r<n$ and}\\
 \qquad\varepsilon_1\!=\!2^{(1+na_n-\frac12b_n)/\sqrt{b_n}}        
             &\text{if $r=n$ } \\
 2^{(1-\frac12a_n)/\sqrt{a_n}}[f_0+\frac{(r+1)f_{i+1}}{a_{n-r-1}}]     
             &\text{if $i=(r+1)a_n-1$}\\
             &\hfil\text{with $r=0,1,\ldots,n-1$}\\
 \varepsilon_1f_{i+1}-b_n\varepsilon_2 f_{i+1-b_n} 
             &\text{ if $i=na_n+rb_n$}\\
 \quad\text{where} 
             &\hfil\text{with $r=1,2,\ldots,n$}\\
 \qquad\varepsilon_2\!=\! 2^{(1+na_n-\frac12b_n)/\sqrt{b_n}}\\
 \qquad\varepsilon_1=2^{(1+na_n-\frac12b_n)/\sqrt{b_n}}           
             &\text{if $r<n$, and} \\
 \qquad\varepsilon_1=2^{(v_n+1-\frac12 a_{n+1})/\sqrt{a_{n+1}}}  
             &\text{if $r=n$}\\
 2^{-((r+1)a_n+\frac12b_n-1)/\sqrt{b_n}}\Bigl[\!\sum_{j=0}^r b_n^j\! f_{i-jb_n+1}
             &\text{if $i=(r+1)(a_n+b_n)-1$}\\
 \hfil+b^{r+1}_n\bigl(f_0+\frac{(r+1)f_{(r+1)a_n}}{a_{n-r-1}}\bigr)\Bigr]
             &\hfil\text{ with $r=0,1,\ldots,n-1$}
        \end{cases}
\end{displaymath}

Inspecting the matrix line by line we observe that, assuming $(a_n)$
and $(b_n)$ are increasing sufficiently rapidly, it follows that
$\norm{S}\leq 2$.  Again by inspecting each line of the matrix, we
deduce that if $f^*_j$ is the $j$-th coordinate functional on
$\ell_1$, $j\geq 0$, it follows that $\lim_{i\to\infty}
f^*_j(S(f_i))=0$. In other words, the rows of the matrix converge to
zero. Therefore $S$ is the adjoint of a linear bounded operator on
$c_0$.

\begin{thm}\label{T:2.1}
  $S$ has a non-trivial closed invariant subspace.
\end{thm}

We shall show that $S$ has an invariant subspace by producing
a vector $x_\infty$ such that the linear span of the orbit of
$x_\infty$ stays away from $e_0$, hence its closure is a non-trivial
$S$-invariant subspace.

We will introduce the following notations.

First we chose  two sequences of positive integers $(m_i)$ and $(r_i)$ as
  follows.  Let $m_0\geq 2$ be arbitrary, put $r_0=1$. Once $m_i$ and $r_i$
  are defined,  choose $r_{i+1}\in\N$ so that
\begin{equation}\label{E:2.2}
r_{i+1}\in[ a_{m_i-1}\cdot\max\limits_{\ell\leq v_{m_i-1}}\norm{\hat e_\ell},
          1+a_{m_i-1}\cdot\max\limits_{\ell\leq v_{m_i-\ell}}\norm{\hat e_\ell}]
\end{equation}
and  let 
\begin{equation}\label{E:2.3}
m_{i+1}=m_i+r_{i+1}.
\end{equation}
Define an increasing sequence $(j_i)$ of positive integers
  inductively: pick any
\begin{equation}\label{E:2.4}  
j_0\in[r_0a_{m_0},r_0a_{m_0}+v_{m_0-r_0}],
\end{equation}
and  once $j_i$ is defined, put
\begin{equation}\label{E:2.5}  
j_{i+1}=j_i+r_ib_{m_i}+r_{i+1}a_{m_{i+1}}.
\end{equation}
Finally, for each $i\geq 0$ define
  \begin{eqnarray}\label{E:2.6}
    p_i &=& \prod\limits_{k=0}^ib_{m_k}^{-r_k};\\
 \label{E:2.7}   z_i &=& f_{j_i+r_ib_{m_i}}+b_{m_i}f_{j_i+(r_i-1)b_{m_i}}+\dots
      +b_{m_i}^{r_i-1}f_{j_i+b_{m_i}}+
      \frac{r_{i+1}f_{j_{i+1}}}{a_{m_i}};\\
 \label{E:2.8}   x_i &=& p_{i-1}\hat e_{j_i}.
  \end{eqnarray}

We note the following easy to prove properties  for our choices.

\begin{prop}\label{P:2.2}
  For each $i\geq 0$ the following statements hold.
  \begin{enumerate}
    \item[a)] $j_i\in[r_ia_{m_i},r_ia_{m_i}+v_{m_i-r_i}]$;
    \item[b)] $x_{i+1}=x_i+p_iz_i$, and thus 
       $x_i=\hat e_{j_0}+\sum_{k=0}^{i-1}p_kz_k$;
    \item[c)]
       If $i$ and $i+\ell$ both belong to $[ra_n,ra_n+v_{n-r}]$ or both belong 
       to $[r(a_n+b_n),na_n+rb_n]$, then $S^\ell f_i=f_{i+\ell}$;
    \item[d)] 
       If $\ell<m_ia_{m_i}-j_i$ then $\min\supp S^\ell z_k\geq j_i+b_{m_i}$ whenever
       $k\geq i$.
  \end{enumerate}
\end{prop}

\begin{proof}
 (a) The proof is by induction. For $i=0$ the required inclusion follows
  from the choice of $j_0$, and if this condition holds for $j_i$ then
\begin{align*}
    j_{i+1}&=j_i+r_ib_{m_i}+r_{i+1}a_{m_{i+1}}\\
            &\in [r_ia_{m_i}+r_ib_{m_i}+r_{i+1}a_{m_{i+1}},
       r_ia_{m_i}+v_{m_i-r_i}+r_ib_{m_i}+r_{i+1}a_{m_{i+1}}]\\
     &\subseteq [r_{i+1}a_{m_{i+1}},r_{i+1}a_{m_{i+1}}+m_i(a_{m_i}+b_{m_i})]=
      [r_{i+1}a_{m_{i+1}},r_{i+1}a_{m_{i+1}}+v_{m_i}].
  \end{align*}

\noindent  
(b) First note that by using ($\widehat{\text{D}}$)  we obtain
 for a $i\in[r(a_n+b_n),na_n+rb_n]$, with $1\leq r\leq n$ in $\N$, that
\begin{align}\label{E:2.2.1}
\hat e_i&=b_n\hat e_{i-b_n}+ f_i\\
        &=b_n^2\hat e_{i-2b_n}+b_nf_{i-b_n}+f_i\notag\\
        &\vdots\notag\\
        &=b_n^r\hat e_{i-rb_n}+b_n^{r-1}f_{i-(r-1)b_n}+\ldots+b_nf_{i-b_n}+f_i.\notag
\end{align}

Note that $j_{i}+r_ib_{m_i}\in[r_i(a_{m_i}+b_{m_i}),m_ia_{m_i}+r_ib_{m_i}]$.
By using first ($\widehat{\text{A}}$) and then (\ref{E:2.2.1}) we obtain
\begin{align*} 
\hat e_{j_{i+1}}&=\hat e_{j_{i}+r_ib_{m_i}+r_{i+1}a_{m_i}}\\
                &=\hat e_{j_{i}+r_ib_{m_i}}+\frac{r_{i+1}}{a_{m_i}} 
                                    f_{j_{i}+r_ib_{m_i}+r_{i+1}a_{m_i}}\\
               &=b_{m_i}^{r_i} \hat e_{j_i}+b_{m_i}^{r_i-1}f_{j_i+b_{m_i}}+\ldots
                   +b_{m_i}f_{j_i+(r_i-1)b_{m_i}}+ 
                          \frac{r_{i+1}}{a_{m_i}}  f_{j_{i}+r_ib_{m_i}+r_{i+1}a_{m_i}}\\
               &=b^{r_i}_{m_i} \hat e_{j_i}+z_i.
\end{align*}
Thus, $x_{i+1}= p_i\hat e_{j_{i+1}}=p_{i-1}\hat e_{j_i}+p_iz_i=x_i+p_iz_i$.

\noindent
(c) If $i$ and $i+\ell$ are both in $[ra_n,ra_n+v_{n-r}]$ it follows
from ($\widehat{\text{A}}$) that
\begin{equation*}
 S^\ell(f_i)=\frac{a_{n-r}}{r} S^\ell(\hat e_i-\hat e_{i-ra_n})=
   \frac{a_{n-r}}{r} (\hat e_{i+\ell}-\hat e_{i-ra_n+\ell})=f_{i+\ell}.
\end{equation*}
The second part of (c) can be deduced in a similar way using  ($\widehat{\text{C}}$). 

\noindent
(d) First note that for $k\geq i$ it follows that (recall that
$m_k\geq m_0\geq 2$)
\begin{equation*} 
  m_ka_{m_k}-j_k>(m_k-r_k-1)a_{m_k}= 
  (m_{k-1}-1)a_{m_k}\geq m_{k-1}a_{m_{k-1}}-j_{k-1}.
\end{equation*}
We can therefore assume that $k=i$. Furthermore, note that for any
$1\leq r\leq r_i$ it follows that
\begin{equation*} 
  r(a_{m_i}+b_{m_i})\leq j_i+r b_{m_i}\leq  j_i+r b_{m_i}+\ell\leq m_i a_{m_i} +r b_{m_i}
\end{equation*} 
and 
\begin{align*} 
  r_{i+1} a_{m_{i+1}}&\leq j_{i+1}\leq j_{i+1}+\ell\leq j_{i+1}+m_i a_{m_i}-j_i\\
    &=r_{i+1} a_{m_{i+1}}+ r_i b_{m_i}+m_i a_{m_i}\leq 
    r_{i+1} a_{m_{i+1}}+v_{m_i}= r_{i+1} a_{m_{i+1}}+v_{m_{i+1}-r_{i+1}}
\end{align*}
Therefore the claim follows from the definition of $z_i$, (\ref{E:2.7})
 and part (c).
\end{proof}
 
Notice that
\begin{equation*}
  \norm{z_i}=1+b_{m_i}+b_{m_i}^2+\dots+b_{m_i}^{r_i-1}+
    \frac{r_{i+1}}{a_{m_i}}\leq
  m_ib_{m_i}^{r_i-1}+\frac{r_{i+1}}{a_{m_i}}
\end{equation*}
Further, since $p_i\leq\frac{1}{b^{r_i}_{m_i}}$, we have
$$\norm{p_iz_i}\leq\frac{m_i}{b_{m_i}}+
\frac{r_{i+1}}{a_{m_i}b_{m_i}^{r_i}}.$$  The
series $\sum_{i=0}^\infty \frac{m_i}{b_{m_i}}$ converges because $(b_i)$
increases sufficiently rapidly.
 Secondly, it follows from the definition of $(r_i)$ that
$$a_{m_i}^{-1}r_{i+1}\leq a_{m_i}^{-1} [1+ a_{m_i-1}
\cdot\max\limits_{\ell\leq v_{m_i-1}}\norm{\hat e_\ell}].$$
Thus,  again since  $(b_i)$ is increasing fast enough, it follows
that
the series
$\sum_{i=0}^\infty\frac{r_{i+1}}{a_{m_i}b_{m_i}^{r_i}}$
converges. Therefore the  $\sum_{i=0}^{\infty}p_iz_i$ converges, and 
the following definition is justified.

\begin{defin} \label{d:x-infty}
  Define $x_\infty=
 \lim_ix_i=\lim_ip_{i-1}\hat e_{j_i}=\hat e_{j_0}+\sum_{i=0}^{\infty}p_iz_i$. 
\end{defin}

Now we can state and prove the key result for proving Theorem~\ref{T:2.1}.

\begin{lem} \label{l:C}
  There exists a constant $C>0$ such that $\dist(y,e_0)\geq C$ for
  every $i$ and every vector of the form
  $y=\sum_{j=j_i}^{m_ia_{m_i}}\gamma_j\hat e_j$.
\end{lem}

\begin{proof}
  Let $C=\inf\Bigl\{\dist(y,e_0)\mid
  y=\sum_{j=j_0}^{m_0a_{m_0}}\gamma_j\hat e_j\Bigr\}$. Since the
  infimum is taken over a finite-dimensional set, it must be
  attained at some~$y_0$. However since all $\hat e_j$ are linear
  independent, it follows that $C=\dist(y_0,e_0)>0$.

  We shall prove the statement of the lemma by induction on $i$. The way
  we defined $C$ guarantees that the base of the induction holds.
  Suppose $y=\sum_{j=j_i}^{m_ia_{m_i}}\gamma_j\hat e_j$. Write
  $y=y_1+y_2+y_3$, where
  \begin{equation*}
    y_1=\!\!\sum_{j=j_i}^{r_ia_{m_i}+v_{m_{i-1}}}{\hskip -.5cm}\gamma_j\hat e_j,\quad
    y_2=\!\sum_{r=r_i+1}^{m_i}
        \sum_{j=ra_{m_i}}^{ra_{m_i}+v_{m_i-r}}{\hskip -.5cm}\gamma_j\hat e_j,
        \quad\text{ and }\quad
    y_3=\!\sum_{r=r_i}^{m_i-1}
        \sum_{j=ra_{m_i}+v_{m_i-r}+1}^{(r+1)a_{m_i}-1}{\hskip -.7cm}\gamma_j\hat e_j.
  \end{equation*}
  Notice that by ($\widehat{\text{B}}$)
  $$y_3=\sum_{r=r_i}^{m_i-1}
        \sum_{j=ra_{m_i}+v_{m_i-r}+1}^{(r+1)a_{m_i}-1}
        \gamma_j 2^{-(r+\frac12-j)/\sqrt{a_{m_i}}}f_j,$$
  so that $\supp y_3\subseteq\bigcup_{r=r_i}^{m_i-1}
   (ra_{m_i}+v_{m_i-r},(r+1)a_{m_i})$.
  Furthermore, using ($\widehat{\text{A}}$), we write $y_2=y_2'+y_2''$ where
  $$y_2'=\sum_{r=r_i+1}^{m_i}
        \sum_{j=ra_{m_i}}^{ra_{m_i}+v_{m_i-r}}\gamma_j\hat e_{j-ra_{m_i}}=
        \sum_{r=r_i+1}^{m_i}
        \sum_{j=0}^{v_{m_i-r}}\gamma_{j+ra_{m_i}}\hat e_j$$ 
  $$\text{and}\quad y_2''=\sum_{r=r_i+1}^{m_i}
          \sum_{j=ra_{m_i}}^{ra_{m_i}+v_{m_i-r}}
          \frac{\gamma_jr}{a_{m_i-r}}f_j.$$
  Therefore,
  $$\supp(y_1+y_2)\subseteq[0,r_ia_{m_i}+v_{m_{i-1}}]
    \cup\bigcup_{r=r_i+1}^{m_i}[ra_{m_i},ra_{m_i}+v_{m_i-r_i}].$$
  One observes that the vectors $y_1+y_2$ and $y_3$ have disjoint
  supports, it follows that
  $\dist(y,e_0)\geq\dist(y_1+y_2,e_0)$.

  Furthermore,
  \begin{equation*}
    \norm{y_2'}=\Bignorm{\sum_{r=r_i+1}^{m_i}
        \sum_{j=ra_{m_i}}^{ra_{m_i}+v_{m_i-r}}\gamma_j\hat e_{j-ra_{m_i}}}\leq
        \sum_{r=r_i+1}^{m_i}
        \sum_{j=ra_{m_i}}^{ra_{m_i}+v_{m_i-r}}\abs{\gamma_j}
        \cdot\max\limits_{k\leq v_{m_{i-1}-1}}\norm{\hat e_k}
  \end{equation*}
  By choice of $(r_i)$  (\ref{E:2.2}), we have
  $\max\limits_{k\leq v_{m_{i-1}-1}}\norm{\hat e_k}\leq
   \frac{r_i}{a_{m_i-r_i-1}}\leq\frac{r}{a_{m_i-r}}$ when
  $r_i<r\leq m_i$. This yields
 \begin{equation*}
    \norm{y_2'}\leq
    \Bignorm{\sum_{r=r_i+1}^{m_i}
        \sum_{j=ra_{m_i}}^{ra_{m_i}+v_{m_i-r}}
        \frac{\gamma_jr}{a_{m_i-r}}f_j}=\norm{y_2''}.
  \end{equation*}
  Since the support of $y_2''$ is disjoint from that of $y_1+y_2'$ and 
  doesn't contain 0, we have
  \begin{align*}
    \dist(y_1,e_0)&\leq\dist(y_1+y_2',e_0)+\norm{y_2'}\\
                  &=\dist(y_1+y_2'+y_2'',e_0)-\norm{y_2''}+\norm{y_2'}\\
                  & \leq\dist(y_1+y_2,e_0) \leq\dist(y,e_0).
  \end{align*}

  It is left to show that $\dist(y_1,e_0)\geq C$. Since $j_i\geq
  r_ia_{m_i}$, it follows from ($\widehat{\text{A}}$) that $y_1=y_1'+y_1''$ where
  $$y_1'=\sum_{j=j_i}^{r_ia_{m_i}+v_{m_{i-1}}}\gamma_j\hat e_{j-r_ia_{m_i}}
    \quad\mbox{ and }\quad
    y_1''=\sum_{j=j_i}^{r_ia_{m_i}+v_{m_{i-1}}}\frac{\gamma_jr}{a_{m_i-r_i}}f_j.$$
  Since $j_i=j_{i-1}+r_{i-1}b_{m_{i-1}}+r_ia_{m_i}$, we have
  $y_1'=\sum_{j=j_{i-1}+r_{i-1}b_{m_{i-1}}}^{v_{m_{i-1}}}\beta_j\hat e_j$,
  where $\beta_j=\gamma_{j+r_ia_{m_i}}$.
   In particular this means, that
  $\supp y_1'\subseteq[0,v_{m_{i-1}}]$, while
  $\min \supp y_1''\geq j_i\geq r_ia_{m_i}$. Thus, the supports are
  disjoint, which yields $\dist(y_1,e_0)\geq\dist(y_1',e_0)$.

  Split the index set of $y_1'$ into two disjoint subsets: let
  \begin{align*}
    A&\!=\![j_{i-1}\!+\!r_{i-1}b_{m_{i-1}},v_{m_{i-1}}]\cap
       \!\bigcup_{r=r_{i-1}}^{m_{i-1}}\!
       \bigl(m_{i-1}a_{m_{i-1}}\!+\!rb_{m_{i-1}},
       (r\!+\!1)(a_{m_{i-1}}\!+\!b_{m_{i-1}})\bigr);\\
    B&\!=\![j_{i-1}\!+\!r_{i-1}b_{m_{i-1}},v_{m_{i-1}}]\cap
       \!\bigcup_{r=r_{i-1}}^{m_{i-1}}\!
       \bigl[r(a_{m_{i-1}}\!+\!b_{m_{i-1}}),m_{i-1}a_{m_{i-1}}\!+\!rb_{m_{i-1}}\bigr].
  \end{align*}
  Write $y_1'=z_a+z_b$ where $z_a=\sum_{j\in A}\beta_j\hat e_j$ and 
  $z_b=\sum_{j\in B}\beta_j\hat e_j$. For $j\in A$ we have
  $\hat e_j=2^{((r+1/2)b_{m_{i-1}}-j)/\sqrt{b_{m_{i-1}}}}f_j$, so that
  $\supp z_a\subseteq A$.
  In view of~(\ref{E:2.2.1}) we can write $z_b=z_b'+z_b''$, where
  $$z_b'=\sum_{j\in B}\sum_{k=0}^{r-1}
      \beta_jb_{m_{i-1}}^kf_{j-kb_{m_{i-1}}}\quad\mbox{ and}\quad
    z_b''=\sum_{j\in B}\beta_jb_{m_{i-1}}^r\hat e_{j-rb_{m_{i-1}}}.$$
  We first note that $\supp z_b'\subseteq B$ and
  determine the support of $z_b''$ as follows.
  If $j\in B$ then $j\geq j_{i-1}+r_{i-1}b_{m_{i-1}}$ and
  $j\in\bigl[r(a_{m_{i-1}}+b_{m_{i-1}}),
   m_{i-1}a_{m_{i-1}}+rb_{m_{i-1}}\bigr]$
  for some $r\in[r_{i-1},m_{i-1}]$. If $r=r_{i-1}$ then
  $j-rb_{m_{i-1}}\geq j_{i-1}$. If $r>r_{i-1}$ then
  $j-rb_{m_{i-1}}\geq ra_{m_{i-1}}>
   r_{i-1}a_{m_{i-1}}+v_{m_{i-2}}\geq j_{i-1}$ by (\ref{E:2.5}).
  We see that $z_b''$ is a linear combination of $\hat e_j$'s with $j_{i-1}\leq j\leq
  m_{i-1}a_{m_{i-1}}$. Hence its support is contained in
  $[0,m_{i-1}a_{m_{i-1}}]$ and, therefore, is disjoint from that of
  $z_a$ and $z_b'$. It follows that
  $\dist(y,e_0)\geq\dist(y_1',e_0)\geq\dist(z_b'',e_0)$. Finally,
  $\dist(z_b'',e_0)\geq C$ by the induction hypothesis.
\end{proof}

\begin{proof}[Proof of Theorem \ref{T:2.1}]
  We will prove that the linear span of the orbit of $x_\infty$ under
  $S$ is at least distance $C$ from $e_0$, hence its closure is a
  non-trivial invariant subspace for $S$.
  Consider a linear combination $\sum_{\ell=0}^N\alpha_\ell S^\ell x_\infty$. It follows
  from (\ref {E:2.5}) that the sequence $(m_ia_{m_i}-j_i)$ is
  unbounded, so that $N<m_ia_{m_i}-j_i$ for some $i\geq 0$. Recall that 
  $x_\infty=x_i+\sum_{k=i}^\infty p_kz_k$, then
  $$\sum_{\ell=0}^N\alpha_\ell S^\ell x_\infty=\sum_{s=0}^N\alpha_\ell S^\ell x_i+
    \sum_{\ell=0}^N\sum_{k=i}^\infty \alpha_\ell S^\ell(p_kz_k).$$
  Notice that the two sums have disjoint supports, and the support of
  the second one doesn't contain $0$. Indeed,
  since $x_i=p_{i-1}\hat e_{j_i}$ then $S^\ell x_i=p_{i-1}\hat e_{j_i+\ell}$ for
  $\ell=1,\dots,N$. Furthermore,
  $$j_i\leq j_i+\ell\leq j_i+N<j_i+(m_ia_{m_i}-j_i)=m_ia_{m_i}.$$
  It follows that $\sum_{\ell=0}^NS^\ell x_i$ is a linear combination of
  $\hat e_j$'s with $j_i\leq j\leq m_ia_{m_i}$. In particular, its
  support is contained in $[0,m_ia_{m_i}]$. On the other hand,
   Proposition \ref{P:2.2} (d) implies that
  $$\min\supp\Bigl(\sum_{\ell=0}^N\sum_{k=i}^\infty S^\ell(p_kz_k)\Bigr)\geq
  j_i+b_{m_i}.$$ 
Therefore, by Lemma~\ref{l:C}
  $$\dist\Bigl(\sum_{\ell=0}^NS^\ell x_\infty,e_0\Bigr)\geq
    \dist\Bigl(\sum_{\ell=0}^NS^\ell x_i,e_0\Bigr)\geq C.$$
\end{proof}

{\footnotesize

\noindent
Department of Mathematics, Texas A$\&$M University, 
College Station, TX 77843, schlump@math.tamu.edu

\noindent 
Department of Mathematics, University of Texas at Austin,
Austin, TX 78712, vladimir@mail.utexas.edu 

}

\end{document}